\documentclass[12pt]{amsart}

\usepackage{amsmath,amssymb,amscd,amsxtra}
\usepackage{latexsym}
\usepackage[dvips]{graphics,epsfig}

\newcommand{\R}{\mathbb{R}}

\newcommand{\Z}{\mathbb{Z}}
\newcommand{\M}[1]{\mathcal{M}^{#1}}
\newcommand{\Mco}[2]{\mathcal{M}_{comp}^{#1}({#2})}
\newcommand{\Mloc}[2]{\mathcal{M}_{loc}^{#1}({#2})}
\renewcommand{\S}{\mathcal{S}}
\newcommand{\di}{\mathcal{S}^{'}}

\newcommand{\nm}[2]{\|{#1}\|_{#2}}

\newcommand{\biggparen}[1]{\biggl(#1\biggr)}

\newcommand{\rd}{{\R}^{d}}

\newcommand{\cF}{\mathcal{F}}
\newcommand{\cFco}[2]{(\mathcal{F} L^{#1})_{comp}({#2})}
\newcommand{\cFloc}[2]{(\mathcal{F} L^{#1})_{loc}({#2})}

\newtheorem{thm}{Theorem}
\newtheorem{lemma}{Lemma}
\newtheorem{prop}{Proposition}

\theoremstyle{remark}
\newtheorem{rem}{Remark}
\theoremstyle{definition}
\newtheorem{deft}{Definition}

\begin{document}

\title{A Beurling-Helson type theorem for modulation spaces}

\author{Kasso A.~Okoudjou}

\address{Kasso A.~Okoudjou\\
Department of Mathematics\\
University of Maryland\\
College Park, MD, 20742 USA}

\email{kasso@math.umd.edu}


\date{\today}


\begin{abstract} We prove a Beurling-Helson type theorem on modulation
spaces. More precisely, we show that the only  $\mathcal{C}^{1}$ changes of
variables that leave invariant the modulation spaces
$\M{p,q}(\rd)$ are affine functions on $\rd$.   A special case of
our result involving the Sj\"ostrand algebra was considered
earlier by A.~Boulkhemair.
\end{abstract}

\maketitle \pagestyle{myheadings} \thispagestyle{plain}
\markboth{K. A. OKOUDJOU}{BEURLING-HELSON TYPE THEOREM}

\section{Introduction}
Given  a function $\phi$ defined from  the torus $\mathbb{T}$ to itself, let
$\phi^{*}$ be the change of variables defined by
\begin{equation}\label{fistar} \phi^{*}(u) = u \circ \phi
\end{equation} for any function $u$ defined on $\mathbb T$.

In $1953$, A.~Beurling and H.~Helson proved that if $\phi$ is
continuous from $\mathbb T$ into itself and if $\phi^{*}$ is a bounded linear operator on  the
Fourier algebra $A(\mathbb {T})= A_{1}(\mathbb {T})$ of absolutely
convergent Fourier series, then necessarily $\phi(t) = kt
+\phi(0)$ for some $k \in \Z$  \cite{BeHe53}. The proof of this
result involved some nontrivial arithmetical considerations. A
different proof was given by J.-P.~Kahane \cite{kaha}. The
Beurling-Helson theorem was later extended to the higher
dimensional setting by W.~M.~Self \cite{self76}. More recently,
V.~Lebedev and A.~Olevski\v{\i} \cite{leole1} further extended and generalized the Beurling-Helson theorem. In particular, for $d\geq 1$ and
$1\leq p < \infty$ let  $A_{p}(\rd)=\mathcal{F} L^{p}(\rd)$
equipped with the norm $\nm{f}{A_{p}}=\nm{\hat{f}}{L^{p}}$ where
$\mathcal F$ is the Fourier transform defined by $\mathcal{F} f
(\omega) = \hat{f}(\omega) = \int_{\rd}f(t)\, e^{-2\pi i t\cdot
\omega}\, dt.$ It was proved in \cite{leole1} that if $\phi : \rd \to \rd$ is $\mathcal C^{1}$, and if
$\phi^{*}$ maps $A_{p}(\rd)$ into itself for some $1\leq p<\infty$,
$p\neq 2$, then $\phi(x) = Ax + \phi(0)$ where $A$ is a real
invertible $d\times d$ matrix. In this higher dimensional setting, the case $p=1$ was already proved in \cite{self76}. Observe that since $A_{2}(\rd)=\mathcal{F}L^{2}(\rd)=L^{2}(\rd)$, the class of functions $\phi$ such that $\phi^{*}$ is bounded on $A_{2}(\rd)$ is quite large. For instance, for any homeomorphism $\phi$ on $\mathbb{T}$ such that $\phi^{-1}$ satisfies the Lipschitz condition, $\phi^{*}$ is bounded on $A_{2}(\mathbb{T})$, and a transference argument can be used to prove similar result for $A_{2}(\rd)$. 

In this note, we shall characterize the $\mathcal{C}^{1}$ changes
of variables that leave invariant the modulation spaces (to be defined below). In particular, our result applies to a special subspace of the
Fourier algebra called Feichtinger algebra. This space denoted
$S_{0}$ was introduced by H.~Feichtinger \cite{Fei81} and is the
smallest Banach algebra that is invariant under both the
translation and the modulation operators. Moreover, the
Feichtinger algebra is an example of a modulation space
 and plays an important role in the theory Gabor
frames \cite{Gr01}. In fact, the modulation spaces have also been
playing an increasing role in the analysis of pseudodifferential
operators \cite{grhe, he03, toft}. Furthermore, a Banach algebra of
pseudodifferential operators known as the Sj\"ostrand algebra,
denoted $S_\omega$, and which contains the H\"ormander class
$S_{0,0}^{0}$, was introduced independently by  Feichtinger \cite{Feiori83} and J.~Sj\"ostrand
\cite{sjo94}. This space is yet
another example of a modulation space.  We refer to \cite{Fei83} for an updated version of \cite{Feiori83} which contains some historical perspectives on the modulation spaces. In $1997$, A.~Boulkhemair \cite{boulk97}
proved that if $\phi$ is a $\mathcal{C}^{1}$ mapping on $\rd$ such
that $\phi^{*}$ maps $S_{\omega}$ into itself, then $\phi$ must be an
affine function: This is a Beurling-Helson type theorem for the Sj\"ostrand algebra. It is therefore natural to seek a
characterization of the changes of variables that leave invariant
modulation spaces. The goal of this note is to extend and
generalize this Beurling-Helson type theorem to all the modulation
spaces. The main argument in the proof of our result is the fact
that the intersection of a modulation space with the space of
functions with compact support coincides with the subspace of
compactly supported functions in $A_{p}(\rd) = \cF L^{p}(\rd)$.
The proof of this fact as well as the definition of the modulation
spaces will be given in Section~\ref{sec2}. Finally in
Section~\ref{sec3} we shall prove our main result. In the sequel, we
shall denote by $|A|$ the Lebesgue measure of a measurable subset
$A$ of $\rd.$

\bigskip

\section{Preliminaries}\label{sec2}

\subsection{Modulation spaces}
The Short-Time Fourier Transform
(STFT) of a function $f$ with respect to a window $g$ is
$$V_g f(x, y)=\int_{\R} f(t)\, \overline{g(t-x)} \, e^{-2\pi i y  t} \, dt,$$
whenever the integral makes sense. This definition can be
extended  to $f\in \S ' (\R^d )$ and $g\in \S (\R^d )$ and yields
a continuous function $V_gf $, see
 \cite{Gr01}.

\begin{deft}
Given $1\leq p, q \leq \infty$, and given a window function $0\neq g \in \S$, the modulation space $\M{p,q}=\M{p,q}(\R^d)$
is the space of all distributions $f \in \di$ for which the following norm is finite:
\begin{equation}
\label{mod}
\nm{f}{\M{p,q}}=\biggparen{\int_{\R^d}\biggparen{\int_{\R^d}|V_{g}f(x,y)|^p\,
dx}^{q/p}\, dy}^{1/q},
\end{equation}
with the usual modifications if $p$ and/or $q$ are infinite.
\end{deft}

\begin{rem}
 The definition is independent of the choice of the window $g$ in
the sense of equivalent norms.

 The modulation spaces were originally introduced by Feichtinger
\cite{Feiori83}.  We refer to \cite{Gr01} and the references therein
for more details about modulation spaces.

 The Feichtinger algebra $S_0$  which coincides with the
modulation space $\M{1, 1}(\rd)$ is  a Banach algebra under both
pointwise multiplication and convolution. Furthermore,
$\M{1,1}(\rd)$ like $\M{p,p}(\rd)$ $1\leq p \leq \infty$ is invariant under the Fourier transform \cite{Fei81, Fei89, Fei90}.

 While the Beurling-Helson theorem completely classifies the
changes of variables that operate in $A_{1}(\mathbb T)$ (and also
on $A_{1}(\R)$) it was still unknown what changes of variables
operate on the Feichtinger algebra and more generally on the
modulation spaces. This question will be completely settled below.
\end{rem}

\subsection{Local modulation spaces}
The theory of modulation can be defined in the general setting of
locally compact Abelian groups \cite{Fei83}. In particular, it can
be shown that for $G=\Z^d$ (or any discrete group), $\M{p,q}(G) =
\ell^{p}(G)$. Similarly, if $G=\mathbb T^{d}$ (or any compact
group), $\M{p,q}(G) = \cF L^{q}(G)$. Here we focus on functions that are locally in a modulation space.

In the sequel we shall denote by $\Mco{p,q}{\rd}$ the subspace of
$\M{p,q}(\rd)$ consisting of compactly supported functions, and by
$\Mloc{p,q}{\rd}$ the space of functions that are locally in
$\M{p,q}(\rd)$. In particular, $u \in \Mloc{p,q}{\rd}$ if and only if for each $g \in \mathcal{C}^{\infty}_{0}(\rd)$ with $supp (g) \subset K$ where $K$ is a compact subset of $\rd,$ we have $u_{K}=g \, u \in \M{p,q}(\rd)$, i.e., $u_{K}
\in \Mco{p,q}{\rd}.$ 

$\cFco{q}{\rd}$ and
 $\cFloc{q}{\rd}$ are defined similarly.

The next result contains the key argument in the proof of our main result. We wish to point out that  some special cases of the result are already known. For instance, the result was proved for $\M{\infty, 1}(\rd)$ in \cite[Theorem
5.1]{boulk97}, while \cite{Fei90} dealt with $\M{p, p}(\rd)$ $1\leq p < \infty$. Furthermore, an independent and different proof of part b.\ of Lemma~\ref{locmodsp}  using convolution relations on generalized amalgam spaces was indicated to us by H.~Feichtinger \cite{Feiper}.

\begin{lemma}\label{locmodsp}
Let $1\leq p, q \leq \infty$. Then the following statements hold

\begin{itemize}
\item[a.] $\Mco{p,q}{\rd} = \cFco{q}{\rd}.$
\item[b.] $\Mloc{p,q}{\rd} = \cFloc{q}{\rd}.$
\end{itemize}
\end{lemma}

\begin{proof}
We shall only prove part a.\ of the result as part b.\ follows from the definition of $\Mloc{p,q}{\rd}$.
Furthermore, to prove a.\ it suffices to show that given a
compact subset $K$ of $\rd$ $\M{p, q}(\rd)|_{K} = \cF L^{q}(\rd)|_{K}$. Note that this last equation holds  not only as set equality, but also as equality of Banach spaces with equivalent norms. 

Let $R>0$ be given and let $u \in \cF L^{q}(\rd)$ such that   $supp (
u) \subset B_{R}(0)$. Let $g \in \mathcal{C}_{c}^{\infty} (\rd)$
with $supp (g) \subset B_{R}(0)$. Then, for each $\omega \in \rd$,
$V_{g}u(\cdot, \omega)$ is supported in $B_{2R}(0)$. Thus, using the fact that $|V_{g}u(x, \omega)| = |V_{\hat{g}}\hat{u}(\omega, - x)|=|\mathcal{F}^{-1}(\hat{u}\cdot T_{\omega}\overline{\hat{g}})(x)|$ we have the following estimates

\begin{align*}
\nm{V_{g}u(\cdot, \omega)}{L^{p}} & \leq |B_{2R}(0)|^{1/p}\nm{V_{g}u(\cdot, \omega)}{L^{\infty}}\\
& = |B_{2R}(0)|^{1/p}\nm{\mathcal{F}^{-1}(\hat{u}\cdot T_{\omega}\overline{\hat{g}})}{L^{\infty}}\\
& \leq  |B_{2R}(0)|^{1/p}\nm{\hat{u}\cdot T_{\omega}\overline{\hat{g}}} {L^{1}}\\
& \leq  |B_{2R}(0)|^{1/p}|\hat{u}|\ast |\hat{g}|(\omega).
\end{align*}
Consequently,
$\nm{V_{g}u}{L^{p, q}(\R^{2d})} \leq |B_{2R}(0)|^{1/p} \nm{\hat {u}}{L^{q}(\rd)} \nm{\hat{g}}{L^{1}(\rd)}, $ that is $$\nm{u}{\M{p, q}(\rd)|_{B_{R}(0)}} \leq C(R, p, q, d)\, \nm{u}{\cF L^{q}(\rd)|_{B_{R}(0)}}.$$
Thus,
$$\cF L^{q}(\rd)|_{B_{R}(0)} \subset \M{p,q}(\rd)|_{B_{R}(0)}.$$ For the converse, let $R>0$ be given and
$u \in \M{p,q}(\rd)$ such that  $supp (u) \subset B_{R}(0).$ 
Let $g \in \mathcal{C}^{\infty}_{0}(\rd)$ such that $g\equiv1$ on $ B_{2R}(0)$. It is trivially seen that for all $x \in B_{R}(0)$
and for all $t \in B_{R}(0)$, $g(t-x)=1$. Thus,
 for all $\omega \in \rd$ and for $x \in B_{R}(0)$,
$$\hat{u}(\omega)\chi_{B_{R}(0)}(x) = \chi_{B_{R}(0)}(x)\, V_{g}u(x, \omega) =\chi_{B_{R}(0)}(x)\,
\int_{B_{R}(0)}u(t)\, e^{-2\pi i t\cdot \omega}\, \overline{g(t-x)}\, dt.$$
Therefore,

$$|B_{R}(0)|^{1/p} |\hat{u}(\omega)| = \nm{\chi_{B_{R}(0)}(\cdot)\, V_{g}u(\cdot, \omega)}{L^{p}}.$$
Hence,
$\nm{\hat{u}}{L^q}\leq |B_{R}(0)|^{-1/p} \nm{V_{g}u}{L^{p, q}},$ that is $$\nm{u}{\cF L^{q}(\rd)|_{B_{R}(0)}} \leq C(R, p, q, d)\, \nm{u}{\M{p,q}(\rd)|_{B_{R}(0)}}.$$  Therefore,
$ \M{p,q}(\rd)|_{B_{R}(0)} \subset \cF L^{q}(\rd)|_{B_{R}(0)}.$ We can now conclude that $$\M{p,q}(\rd)|_{B_{R}(0)} = \cF L^{q}(\rd)|_{B_{R}(0)}.$$
\end{proof}

\section{Main results}\label{sec3}
Before stating our main result, we wish to indicate that it is trivially seen that all the modulation spaces are invariant under affine changes of variables. That is, let $1\leq p, q \leq \infty$ and $\phi : \rd \to \rd$ be an affine
mapping, i.e., $\phi(x) = Ax +b$ where $A$ is a $d\times d$ real
invertible matrix and $b \in \rd$. Then the linear operator
$\phi^{*}$ given by~\eqref{fistar} maps $\M{p,q}(\rd)$ into
itself, that is
$$\phi^{*}(\M{p, q}(\rd)) \subset \M{p, q}(\rd).$$
Indeed, let $g\in \S$ and $u \in \M{p,q}(\rd)$, and  $\tilde{g} = g\circ
A^{-1}$ where $A^{-1}$ is the inverse of $A$. The result follows
from
$$V_{g}\phi^{*}(u) (x, \omega) = \tfrac{1}{|det A|} e^{-2\pi i \omega \cdot A^{-1} b}\,
V_{\tilde{g}}u(Ax + b, (A^{*})^{-1} \omega)$$ where $A^{*}$ denote
the conjugate of $A$.

If we restrict our attention to the modulation spaces $\M{p, p}(\rd)$ $1<p<
\infty$ the following stronger result can be proved.
For Proposition~\ref{pwaffmpq} we assume that $\rd = \cup_{k=1}^{N} Q_{k}$
where for each $k,$ $Q_{k}$ is a (possible infinite) ``cube'' with
sides parallel to the coordinates axis. Moreover, we assume that
for $k=1, \hdots, N$ the $Q_k$s have disjoint interiors.

\begin{prop}\label{pwaffmpq}
 Let $\phi$ be a continuous on $\R^d$ such that for
$k=1, \hdots , N$, the restriction $\phi_k$ of $\phi$ to $Q_k$ is
an affine function given by $\phi_k (x) = A_{k} x +b_k$ where
$A_k$ is a real invertible $d\times d$ matrix and $b_k \in \rd$.
Then the linear operator $\phi^{*}$ given by~\eqref{fistar} maps
$\M{p,p}(\rd)$ into itself, that is
$$\phi^{*}(\M{p, p}(\rd)) \subset \M{p, p}(\rd).$$
\end{prop}

\begin{proof}
It is evident from the definition of the modulation spaces that
$\M{p,p}$ is invariant under the Fourier transform, see \cite{Fei89, Fei90}. Let $u \in
\M{p,p}(\rd)$, then $$\phi^{*}(u) = u \circ \phi =\sum_{k=1}^{N}
\chi_{Q_{k}} \cdot (u \circ \phi)=\sum_{k=1}^{N} \chi_{Q_{k}}
\cdot (u \circ \phi_{k}),$$ and so $$\nm{\phi^{*}(u)}{\M{p,p}}
\leq \sum_{k=1}^{N} \nm{\chi_{Q_{k}} \cdot (u \circ
\phi_{k})}{\M{p,p}}.$$ As indicated above, 
$u\circ \phi_k \in \M{p,p}$. Hence, $v_{k}=\mathcal{F}^{-1}(u\circ \phi_{k}) \in \M{p,p}(\rd)$ as well.  Moreover, note that  $\chi_{Q_{k}}$ is a
bounded Fourier multiplier on all $\M{p,p}(\rd)$: this follows
from \cite[Theorem 1]{bggo} in the case $d=1$, and from
\cite[Theorem 6]{bgor} when $d>1$. Consequently, using the invariance of $\M{p,p}(\rd)$ under the Fourier transform, we conclude that there exists $c_k
>0$ such that
$$\nm{\chi_{Q_{k}}\cdot (u \circ \phi_{k})}{\M{p,p}} =
\nm{\mathcal{F}^{-1}(\chi_{Q_{k}}\cdot \hat{v}_{k})}{\M{p,p}} \leq c_k \, \nm{u}{\M{p,p}},$$  from which
the proof follows.
\end{proof}

\medskip

\begin{rem}
The conclusion of Proposition~\ref{pwaffmpq} holds if we used an infinite decomposition of $\rd$, that is if we assume that $\rd= \cup_{k=1}^{\infty} Q_{k}$ where the cubes $Q_{k}$ still have sides parallel to the coordinate axis and disjoint interiors. In this case, the extra assumption needed to prove the previous result is that the constants $c_k$ appearing in the above proof, are uniformly bounded, i.e., $\sup_{k}c_{k} < \infty.$

\end{rem}

\medskip

We are now ready to state and prove our main result.

\begin{thm}\label{mpqconv}
Let $\phi : \R^d \to \R^d$ be a $\mathcal{C}^{1}$ function. Assume that the
operator $\phi^{*}$ defined by~\eqref{fistar} maps $\M{p,q}(\rd)$
into itself, i.e., $\phi^{*}(\M{p,q}(\R^d)) \subset \M{p,q}(\R^d)$
for some $1\leq p, q \leq \infty$, with $2\neq q < \infty$. Then
$\phi$ is an affine mapping, that is $\phi(x) = Ax + \phi(0)$ for
some real invertible $d\times d$ matrix $A $.

In particular, the Feichtinger algebra $\M{1,1}(\R^d)$ is
preserved by, and only by affine changes of variables.
\end{thm}

\begin{proof}
Because $\phi^{*}(\M{p,q}(\R^d)) \subset \M{p, q}(\R^d)$ and
$\phi^{*}(u) = u \circ \phi$ is compactly supported whenever $u$
is, Lemma~\ref{locmodsp} implies that $\phi^{*}$ maps
$\Mco{p,q}{\R^d}=\cFco{q}{\R^d}$ into itself as well as
$\Mloc{p,q}{\R^d}=\cFloc{q}{\R^d}$ into itself. Therefore,

when $d=1$ and $q=1$, the Beurling-Helson Theorem \cite[pp.
84-86]{BeHe53}, implies  that $\phi(x) = ax + \phi(0)$;

when $d=1$ and $1<q<\infty$, $q\neq 2$, it follows from
\cite[Theorem 3]{leole1} that $\phi(x) = ax + \phi(0)$;

when $d>1$ and $q=1$, it follows from \cite[Corollary 1]{self76}
that $\phi(x)= Ax + \phi(0)$, where $A$ is a real invertible
$d\times d$  matrix;

when $d>1$ and $1< q < \infty$, $q\neq 2$, it follows from
\cite[Theorem 6]{leole1} that $\phi(x) = Ax + \phi(0)$,
 where $A$ is a real invertible $d\times d$  matrix.
\end{proof}

\begin{rem}
The fact that $q\neq2$ in Theorem~\ref{mpqconv} was justified in the Introduction. Moreover, we restricted to $q<\infty$, because the key ingredients in the proof of our main result are \cite[Theorem 3, Theorem 6]{leole1} whose proofs are based on a density argument. It is not clear to us if Theorem~\ref{mpqconv} holds for $q=\infty$. 
\end{rem}

\begin{rem}
Using Lemma~\ref{locmodsp} and \cite[pp.~214]{leole1}, it
follows that if $\phi: \rd \to \rd$ is nonlinear and
$\mathcal{C}^2$, then $\phi^{*}$ is not bounded on $\M{p,q}$. This
fact together with Proposition~\ref{pwaffmpq}, show that the
$\mathcal{C}^1$ condition in Theorem~\ref{mpqconv} is the only
nontrivial smoothness condition to impose on $\phi$.

For the Sj\"ostrand algebra $S_{\omega}$ which coincides with the modulation space $\M{\infty, 1}(\rd)$,  Theorem~\ref{mpqconv}  was proved in under  a weaker assumption on $\phi$.
More specifically, it was proved in \cite[Theorem 5.1]{boulk97} that if $\phi$ is a proper
mapping, i.e., $\phi$ is continuous on $\rd$ and $\phi^{-1}(K)$ is a compact set for any compact subset $K$ of $\rd$,
and if $\phi^{*} (\M{\infty, 1}(\rd)) \subset \M{\infty, 1}(\rd)$ then $\phi(x) = Ax +\phi(0)$.
It is also straightforward to prove Theorem~\ref{mpqconv} under this weaker assumption on $\phi$.

Finally, we wish to conclude this paper by pointing out the connection of our main result to certain Fourier multipliers. More precisely, let  $\sigma$ be a function defined on $\rd$. The Fourier
multiplier with symbol $\sigma$ is  the   operator $H_{\sigma}$
initially defined on $\S$ by $$H_{\sigma}f(x) =
\int_{\rd}\sigma(\xi)\, \hat{f}(\xi) \, e^{2\pi i \xi \cdot x}\,
d\xi.$$ We refer to \cite{Ste70} for  more on Fourier multipliers.
As mentioned above, there is a strong connection between the $L^p$-continuity of the
Fourier multipliers and the Beurling-Helson theorem. In particular, the family of
homomorphisms $e^{i \phi(\xi)}$ on the space of $L^p$-Fourier
multipliers was investigated by H\"ormander in \cite[Section 1.3]{Ho63}.  It is easily seen that
$\sigma_{0}(\xi)
  = e^{i \xi}$, then $H_{\sigma_{0}}$ is bounded on all $L^{p}(\rd)$
  for $1\leq p \leq \infty$ and $d\geq 1$. H\"ormander proved
  that if $\phi: \rd \to \rd$ is $\mathcal{C}^2$ and if
$\phi^{*}(\sigma_{0})(\xi)=\sigma_{0}(\phi (\xi)) = e^{i\phi(\xi)}$ is a bounded Fourier multiplier
  on  $L^{p}(\rd)$ for some $1<p<\infty$ and $p\neq 2$, then
  $\phi$ is an affine function \cite[Theorem 1.15]{Ho63}. It is
  interesting to note that there exist nonlinear (non-affine) functions $\phi$ on $\rd$ such that
the Fourier multipliers with symbols $\phi^{*}(\sigma_{0})(\xi) =\sigma_{0}(\phi(\xi))= e^{i\phi(\xi)}$ are bounded on all modulation spaces \cite{bgor}.
\end{rem}

\bigskip

\section{Acknowledgment} The author would like to thank Chris Heil for bringing some of the questions  discussed in this work to his attention. He also thanks \'Arp\'ad B\'enyi, Hans Feichtinger, Karlheinz Gr\"ochenig, Norbert Kaiblinger, and Luke Rogers for helpful discussions.

\end{document}